\documentclass[a4paper,12pt,twoside]{article}

\usepackage[T1]{fontenc}
\usepackage[utf8]{inputenc}
\usepackage[english]{babel}
\usepackage[parfill]{parskip}

\usepackage{amssymb}
\usepackage{amsmath}
\usepackage{amsthm}
\usepackage{mathtools}
\usepackage[cal=boondox]{mathalfa}
\usepackage{mathrsfs}
\usepackage{dsfont}

\usepackage{graphicx}
\usepackage{caption}
\usepackage{hyperref}
\usepackage{tikz}
\usepackage{pgfplots}
\pgfplotsset{compat=1.18}

\usepackage[bottom=2cm, right=2cm, left=2cm, top=2cm]{geometry}
\usepackage{fancyhdr}

\fancypagestyle{plain}
{
	\fancyhf{}
	
	\fancyhead[CO]{\textsc{Kobayashi distance to submanifolds in a Riemannian manifold}}
	\fancyhead[CE]{\textsc{K. Poirot}}
	\fancyfoot[C]{\thepage}
}

\fancypagestyle{firstpage}[plain]
{
	\fancyhead[C]{}
}

\allowdisplaybreaks % Do not remove, otherwise layout issues.

\usepackage{setspace}
\usepackage{titlesec}

\titleformat{\section}
{\filcenter}{\textbf{\textsc{\Roman{section}.}}}
{.5em}{\bfseries\textsc}

\titleformat{\subsection}
{\filcenter}{\textbf{\textsc{\Roman{section}.\arabic{subsection}.}}}
{.5em}{\bfseries\textsc}

\newtheoremstyle{custom}{}{}{\itshape}{}{\bfseries}{:}{\newline}{}
\newtheoremstyle{custom2}{}{}{}{}{\itshape}{.}{.5em}{}

\theoremstyle{custom}
\newtheorem{Def}{Definition}[section]
\newtheorem*{Def*}{Definition}
\newtheorem{Thm}[Def]{Theorem}
\newtheorem*{Thm*}{Theorem}
\newtheorem{Pro}[Def]{Proposition}
\newtheorem*{Pro*}{Proposition}
\newtheorem{Lem}[Def]{Lemma}
\newtheorem*{Lem*}{Lemma}
\newtheorem{Cor}[Def]{Corollary}
\newtheorem*{Cor*}{Corollary}
\newtheorem*{Cjt*}{Conjecture}
\newtheorem{THM}{Theorem} %For global theorem, no connected to counter

\theoremstyle{custom2}
\newtheorem{Rmk}[Def]{Remark}
\newtheorem{Ex}[Def]{Example}
\newtheorem{Not}[Def]{Notation}

\newcommand{\parag}[1]{\hspace{-\parindent}\textit{#1}.\hspace{.5em}}

\newcommand{\RR}{\mathds{R}}
\newcommand{\CC}{\mathds{C}}
\newcommand{\DD}{\mathds{D}}
\newcommand{\BB}{\mathds{B}}
\newcommand{\HH}{\mathds{H}}
\let\SS\relax\newcommand{\SS}{\mathds{S}}

\newcommand{\Mc}{\mathcal{M}}
\newcommand{\Cc}{\mathcal{C}}
\newcommand{\Sc}{\mathscr{S}}
\newcommand{\Oc}{\mathscr{O}}
\newcommand{\Uc}{\mathcal{U}}
\newcommand{\Vc}{\mathcal{V}}

\DeclareMathOperator{\Hess}{Hess}
\DeclareMathOperator*{\Span}{Span}
\DeclareMathOperator{\tr}{tr}
\DeclareMathOperator{\codim}{codim}

\begin{document}

\begin{center}
\large\textbf{\textsc{Kobayashi distance to submanifolds in a Riemannian manifold}}\normalsize\\
\vspace{0.5cm}
\textsc{Keeran Poirot}
\end{center}
\thispagestyle{firstpage}

\vspace{1cm}

\begin{abstract}
	We prove that totally geodesic hypersurfaces are at infinite Kobayashi distance in a Kobayashi hyperbolic Riemannian manifold. To do this, we prove a Schwarz lemma on minimal maps with a uniform control on their Laplacians. Finally we give a full characterisation of submanifolds at infinite distance on spaces of constant sectional curvature in dimension at least 4.
\end{abstract}

\section*{Introduction}

\hspace{\parindent} The Kobayashi-Royden pseudometric is a useful tool to study the geometry of (almost-)complex manifolds. In recent years there has been a growing interest in possible extensions of such objects to other fields in geometry. It was first extended by Drnov\v{s}ek-Forstneri\v{c} in \cite{dr21} on Euclidean domains, then by Gaussier-Sukhov in \cite{hg23} and \cite{hg24} on arbitrary Riemannian manifolds. More recently Broder-Iliashenko-Madnick in \cite{br25} also extended the Kobayashi-Royden pseudometric on calibrated manifolds. In this paper we will focus on Riemannian manifolds. A natural question arises from such extension: what informations on the geometry of a Riemannian manifold does the Kobayashi metric provide ?

The Kobayashi-Royden pseudometric at a point $p$ and in a direction $\xi$ is obtained by taking the inverse of the maximal radius of holomorphic disks $u$ such that $u(0)=p$ and $du_0(1)=\xi$. This pseudometric can be better understood as a tool which encodes the maximal size of a holomorphic disk to fit in a manifold. Using Schwarz-type lemmas \cite{yau78} on negatively curved complex manifolds, one can control holomorphic disks and show that the pseudometric is non-degenerate and the integrated pseudodistance is a distance. Similar lemmas arise naturally when studying negatively curved Riemannian manifolds. Such manifolds are called Kobayashi hyperbolic, or just hyperbolic if there is no possible confusion. Moreover, a complex manifold is said to be complete hyperbolic if the integrated pseudodistance is a complete distance.

In this paper we will look more precisely at what properties of submanifolds can be removed while preserving the completeness, on a complete hyperbolic Riemannian manifold. This type of question allows us to better understand what properties of submanifolds is encoded with Kobayashi's notion of hyperbolicity. Note that this question has been answered in complex situations in \cite{iv04}. Ivashkovich-Rosay studied submanifolds which can be removed on a Kobayashi hyperbolic complex manifolds while being \emph{detectable} by the completness (near the submanifold) of the Kobayashi-Royden metric on the complementary. They showed that a real submanifold is at infinite Kobayashi distance in a Kobayashi hyperbolic complex manifold if and only if it is either a real Levi-flat hypersurface or a complex hypersurface. Precisely, they showed this result on almost-complex manifolds; however when the complex structure is not integrable, the converse is lost: a hypersurface at infinite distance might not be Levi-flat. Note that, for instance, understanding the behaviour of submanifolds relative to the Kobayashi-Royden metric could be useful when looking at Green's theorem which removes multiple submanifolds on a compact complex manifold and ask if the obtained manifold is complete hyperbolic. See \cite{gr72} and \cite{gr77}. Notably, the key trick to show Green's theorem is to work locally on small holomorphic polydisks around the removed submanifolds which give small complete hyperbolic domains.

The Kobayashi-Royden pseudometric is linked to the idea of a Schwarz lemma: if in a given context there exists some Schwarz-like lemma on a set of disks, one could hope for a similar metric to exist (with, of course, some adequate regularity and existence properties for the set of disks). In particular, on negatively curved Riemannian manifolds, there exists such lemmas controlling minimal disks similarly to the ones in complex geometry. See \cite{fo24} for Euclidean domains, \cite{jo84} and \cite{ch84} for arbitrary Riemannian manifolds. This motivates the study of analogous pseudometrics in a Riemannian context, an idea that has already been proposed by Gromov in \cite{gr07}.

As mentioned earlier, the Kobayashi-Royden pseudometric was first extended on Euclidean domains by \cite{dr21} and later on arbitrary Riemannian manifold by Gaussier-Sukhov in \cite{hg23} who showed some regularity property for the pseudometric which are key to be able to talk about Kobayashi pseudodistance and Kobayashi hyperbolicity. The definition of the Kobayashi-Royden pseudometric on a Riemannian manifold is extremely similar to the complex pseudometric. Instead of being based on holomorphic disks, we use conformal harmonic disks. Then we can define all related notions: Kobayashi hyperbolicity, complete hyperbolicity, ... . Note that the first result in \cite{iv04} was extended in \cite{hg25} where Gaussier-Sukhov showed that a strictly pseudo-convex relatively compact domain on a Riemannian manifold is complete hyperbolic.

\parag{Results} Let us now state the main results of this paper.

\begin{THM}\label{thm:1}
	\hspace{\parindent} Let $(\Mc,g)$ be a Kobayashi hyperbolic Riemannian manifold, $\Sigma$ a hypersurface and $p\in\Sigma$. If $\Sigma$ is locally totally geodesic around $p$, then $p$ is at infinite Kobayashi distance.	
\end{THM}

Here locally totally geodesic around $p$ means that there exists $\Uc\subseteq\Sigma$ a neighborhood of $p$ with vanishing second fundamental form. Also $p$ is at infinite Kobayashi distance if it is not at finite distance; that is, if there exist a fixed point $q$ and a sequence $p_k\rightarrow p$ in the complementary such that the distances $d_{\Mc-\Sigma}(q,p_k)$ are bounded.

Note that the Euclidean situation was already studied by Drnov\v{s}ec-Forstneri\v{c} in \cite{dr21} who showed (Lemma 5.2) that a hyperplane is at infinite distance for Euclidean domains with minimal metric.

Moreover, Theorem \ref{thm:1} is a local result and can be rephrased as: if $\Sigma$ is totally geodesic then every point of $\Sigma$ is at infinite distance from the complementary. Therefore, from a global viewpoint, the following corollary about completeness follows.

\begin{Cor*}
	\hspace{\parindent} Let $(\Mc,g)$ be a complete hyperbolic Riemannian manifold and $\Sigma$ a hypersurface of $\Mc$. If $\Sigma$ is totally geodesic then $\Mc-\Sigma$ is complete hyperbolic.
\end{Cor*}

We then prove that the condition of being locally totally geodesic is also a necessary condition, provided $(\Mc,g)$ has constant sectional curvature.

\begin{THM}\label{thm:2}
	\hspace{\parindent} Let $(\Mc^n,g)$ be a Kobayashi hyperbolic Riemannian manifold with dimension $n>3$ and constant sectional curvature. Let $\Sigma$ be a submanifold of $\Mc$ and $p\in\Sigma$. If $\Sigma$ is not, near $p$, a totally geodesic hypersurface then there is a point $p'\in\Sigma$ near $p$ which is at finite Kobayashi distance from the complementary.
\end{THM}

We will prove Theorem \ref{thm:2} on hyperbolic domains of a model space $\Mc_\kappa^n$, the unique complete simply connected Riemannian $n$-manifold of constant sectional curvature $\kappa$, of dimension strictly greater than $3$. Indeed, the arguments are local and every Riemannian manifold with constant sectional curvature is locally isometric to a model space.

Note that when we say that there exists a point $p'$ near $p$, we mean that there exists a point $p'$ in any open neighborhood of $p$. Again, we have the following corollary from Theorem \ref{thm:1} and Theorem \ref{thm:2} in terms of completness.

\begin{Cor*}
	\hspace{\parindent} Let $(\Mc,g)$ be a complete hyperbolic Riemannian manifold with dimension $n>3$ and constant sectional curvature. Let $\Sigma$ be a submanifold of $\Mc$. Then $\Mc-\Sigma$ is complete hyperbolic if and only if $\Sigma$ is a totally geodesic hypersurface.
\end{Cor*}

\parag{Discussion} Let us outline the proof of both theorems:
\begin{itemize}
	\item[--] For Theorem \ref{thm:1}, we will first show that the signed distance function from a minimal disk is an almost harmonic map. Using a Harnack inequality, we prove a Schwarz lemma, which will allow us to control the behaviour of minimal disks approaching $\Sigma$.
	\item[--] For Theorem \ref{thm:2}, we will first reduce the situation to study \emph{non-minimal} surfaces in $\Mc^3_\kappa$ and then construct adequate sequences of minimal disks approaching the surface.
\end{itemize}
At first, Theorem \ref{thm:2} can be surprising because, from the complex results, one could expect that minimal submanifolds are good candidates to be at infinite Kobayashi distance. Indeed, in the complex settings completeness detects complex properties. Therefore, one could expect to detect some minimality properties in a Riemannian manifold. At the end of this paper, from the proof of Theorem \ref{thm:2} we will give a class of minimal counter exemples.

The existence of totally geodesic hypersurfaces is a restrictive condition, even locally. Actually, on arbitrary Riemannian manifolds such hypersurfaces typically do not exist. This is another surprising result compared to complex geometry. Note that a lot of manifolds still have totally geodesic hypersurfaces, specifically those with local symmetries. Therefore Theorem \ref{thm:1} can be used on such manifolds and not only on the usual model spaces $\Mc_\kappa^n$.

\parag{GenAI disclosure} No generative AI tools were used in the preparation and writing of the proofs and manuscript. Such tools were used to fix an issue on the layout with latex.

\parag{Organization} In Section \ref{sec:1}, we provide some background on the Kobayashi-Royden pseudometric and Kobayashi hyperbolic spaces. In Section \ref{sec:2}, we will look into the proof of Theorem \ref{thm:1} using Fermi's coordinates and a Schwarz-type lemma. Finally, in Section \ref{sec:3}, we prove Theorem \ref{thm:2} and give examples of minimal hypersurfaces at finite distance.

%\renewcommand{\contentsname}{}
%\vspace{-1.3cm}
%\tableofcontents

\parag{Acknowledgements} I would like to thank Hervé Gaussier and Andrea Seppi for drawing my attention to this problem and also for their valuable guidance.

\section{Context and definitions}
\label{sec:1}

\hspace{\parindent} Let $\Mc$ be a Riemannian manifold with metric $g$ and $\Sigma$ be a submanifold of $\Mc$. Recall that a smooth map $f:(\Sc,h)\rightarrow(\Mc,g)$ between two Riemannian manifolds is:
\begin{itemize}
	\item[\textbf{i)}] \textit{conformal} if there exists a smooth function $\lambda:\Sc\rightarrow\RR_{>0}$ such that $f^*g=\lambda^2h$. In particular a conformal map is always an immersion (with this definition).
	\item[\textbf{ii)}] \textit{harmonic} if its tension field $\tau(f):=\tr_h\nabla df$ vanishes.
\end{itemize}
In the special case where $\Sc$ is a Riemann surface both definitions do not depend on the choice of conformal metric on $\Sc$. Given a conformal map $u:\Sc\rightarrow\Mc$, $u$ is harmonic if and only if its image $u(\Sc)$ is an immersed minimal surface, in the usual sense where the trace of the second fundamental tensor $\mbox{I\!I}$ vanishes. Hence, a conformal harmonic map $\DD\rightarrow\Mc$ will be called a \textit{minimal disk}, and the set of minimal disks is denoted by $\mathcal{m}(\Mc)$. For more details on those notions see \cite{jo05}.

\begin{Def}[Kobayashi-Royden pseudometric]
	\hspace{\parindent}The Kobayashi-Royden pseudometric of $\Mc$ is defined by:
	$$\begin{array}{rcl}
	k_\Mc:T\Mc &\longrightarrow& \RR_+ \\[1em]
	(p,v) &\longmapsto& \inf\left\{\dfrac{1}{r},u\in\mathcal{m}(\Mc)\mbox{ such that }u(0)=p,du_0(1)=rv\right\}
	\end{array}.$$
\end{Def}

As mentioned in the introduction, the idea of a pseudometric from minimal disks was first mentioned by Gromov in \cite{gr07}. However the key property of the pseudometric being integrable was proven by Gaussier-Sukhov in \cite{hg23}.

\begin{Thm}[Gaussier-Sukhov \cite{hg23}]
	\hspace{\parindent}There exists a Kobayashi-Royden pseudometric on Riemannian manifolds and it is an upper semicontinuous function on the tangent bundle $T\Mc$.
\end{Thm}

Being upper semicontinuous, the pseudometric can be integrated on paths. In particular taking the infimum of the length of all paths connecting two points gives us a pseudodistance that is a symmetric non-negative function on the set of pair of points vanishing on the diagonal and verifying a triangle inequality. In other words, a pseudodistance is a distance without the separation property. In particular a pseudodistance can vanish where a distance would not.
	
\begin{Def}[Notions of hyperbolicity]
	\hspace{\parindent} We say that $\Mc$ is:
	\begin{itemize}
		\item[--] Royden hyperbolic if for every point $p\in\Mc$ there is an open neighborhood $\Uc$ of $p$ in $\Mc$ and a positive constant $c$ such that for every $(q,v)\in T\Uc,k_\Mc(q,v)\geq c\vert v\vert_{g_q}$.
		\item[--] Kobayashi hyperbolic if the integrated pseudodistance $d_\Mc$ is a distance.
	\end{itemize}
\end{Def}

Moreover, just as in the complex theory, we have the following theorem, proved in \cite{hg23}.
\begin{Thm}[Royden's theorem]
	\hspace{\parindent}$\Mc$ is Royden hyperbolic if and only if it is Kobayashi hyperbolic.
\end{Thm}

Finally, we have one last definition.
\begin{Def}[Complete hyperbolicity]
	\hspace{\parindent}We say that $\Mc$ is complete hyperbolic if $(\Mc,d_{\Mc})$ is complete as a metric space.
\end{Def}

\begin{Ex}
	In \cite{fo24} it is shown that the Euclidean ball is Kobayashi hyperbolic and that the pseudometric coincides with the Beltrami-Cayley-Klein metric. Note that the pseudometric in \cite{fo24} is defined with weakly conformal harmonic disks. But in the case of the Euclidean ball, the extremal disks (disks that realise the infimum) are immersed so we have equality with the pseudometric used in this paper. The real hyperbolic space $\HH^n$ is also Kobayashi hyperbolic with pseudometric equal to its Riemannian metric up to a factor. In particular, those exemples show that the Kobayashi pseudometric is a Riemannian object and not a conformal one, as the hyperbolic and Euclidean metrics are conformal on the ball $\BB^n$ but the Kobayashi-Royden pseudometrics are not.
\end{Ex}

\begin{Rmk}
	Note that the pseudodistance can be better understood using S. Kobayashi's initial construction with chains of disks. See \cite{dr21} and \cite{hg23}. Given two points $p,q\in\Mc$, a \textit{chain of disks} $\alpha$ on $\Mc$ linking $p$ to $q$ is a sequence of minimal disks $u_1,...,u_n:\DD\rightarrow\Mc$, points $p_0=p,p_1,...,p_n=q\in\Mc$ and $a_1,...,a_n,b_1,...,b_n\in\DD$ such that $\forall i\in[\![1,n]\!],u_i(a_i)=p_{i-1}$ and $u_i(b_i)=p_i$. One can also define the length of a chain of disks $\alpha$ as $\ell(\alpha):=\sum_i\rho_\DD(a_i,b_i)$ where $\rho_\DD$ is half of the usual Poincaré distance on $\DD$. Finally, the integrated pseudodistance between two points $p,q\in\Mc$ can be obtained by the following property $d_{\Mc}(p,q)=\inf\ell(\alpha)$ where the infimum is taken over all chain of disks $\alpha$ linking $p$ to $q$. This property will be particularly useful to find counterexemples later.
\end{Rmk}

Let us give some important properties of the Kobayashi pseudometric that we will often use in the proof of Theorems \ref{thm:1} and \ref{thm:2}.
\begin{Pro}
	\hspace{\parindent}Let $\Mc$ and $\Sc$ be two Riemannian manifolds.
	\begin{itemize}
	\item[--] Given a totally geodesic immersion $\iota:\Sc\rightarrow\Mc$, the Kobayashi pseudometric is decreasing along the immersion. Specifically, for all $(p,v)\in T\Sc,k_\Mc(\iota(p),d\iota_pv)\leq k_\Sc(p,v)$. In particular, if $\Sc=\Omega$ is a domain of $\Mc$ then $k_\Mc\leq k_\Omega$. Therefore if $\Mc$ is Kobayashi hyperbolic then $\Omega$ is also Kobayashi hyperbolic
	\item[--] Given $u\in\mathcal{m}(\Mc)$ we have the decreasing property, that is for every $(z,\xi)\in\DD\times\CC$,$k_\Mc(u(z),du_z(\xi))\leq\rho_\DD(z,\xi)$, where $\rho_\DD$ is half the Poincaré metric, that is $\rho_\DD(z,\xi):=\dfrac{\vert\xi\vert}{1-\vert z\vert^2}$.
\end{itemize}
\end{Pro}

\begin{proof}
	Both properties are easy to show using the definitions and basic properties on minimal disks. The first property comes from the fact that minimal disks are preserved under the composition of a totally geodesic immersion. See \cite{ee64}. For the second property, the trick is to use an automorphism of the disk to reparametrize the minimal disks $u$ into a minimal disks $\tilde{u}$ such that $u(z)=\tilde{u}(0)$. From there, we just have to apply the definition of the Kobayashi-Royden pseudometric.
\end{proof}

Now, $\Mc$ will be supposed Kobayashi hyperbolic so $d_{\Mc}$ and $d_{\Mc-\Sigma}$ are distances. Similarly to \cite{iv04}, we have the following definition:

\begin{Def}[Point at infinite distance]
	\hspace{\parindent}A point $p\in\Sigma$ is said to be at finite distance to $\Mc-\Sigma$ if there is a sequence of points $p_k\rightarrow p,p_k\in\Mc-\Sigma$, a fixed point $q\in\Mc-\Sigma$ such that the sequence $(d_{\Mc-\Sigma}(p_k,q))_k$ is bounded. Otherwise, $p$ is said to be at infinite distance.
\end{Def}

\begin{Rmk}\label{rmk:infinite}
	There exists a Kobayashi hyperbolic manifold $\Mc$, a hypersurface $\Sigma$, a point $p\in\Sigma$ at finite distance and a sequence $p_k\in\Sigma\rightarrow p$ such that $(d_{\Mc-\Sigma}(p_k,q))_k$ is not bounded. The simplest example is obtained by taking a strictly convex surface in a Kobayashi hyperbolic domain in $\RR^3$. Note that this example could be done in higher dimension with similar arguments. From one size, one can use a result from Gaussier-Sukhov \cite{hg25} and show that $\Omega$ is complete hyperbolic because strict convexity implies strict minimal pseudoconvexity. But from the other side, we will see in Section \ref{sec:3}, that there is a point $p$ of the surface at finite distance.
\end{Rmk}

\section{Proof of Theorem 1}
\label{sec:2}

\hspace{\parindent} Let us give a quick idea of the proof of Theorem \ref{thm:1}, which follows the proof structure of \cite{iv04}. First, we can reduce the situation by looking at a well-chosen coordinates neighborhood in $\Mc$ of the point $p\in\Sigma$ of interest. While, in general, components of minimal disks will not be harmonic, in those special coordinates we will show that when the disks approach $\Sigma$ one of their components will become more and more harmonic as long as $\Sigma$ is totally geodesic. To be more precise, we will have a control on the laplacian of such components that goes to $0$ uniformly as the disks approach $\Sigma$. Then by using some standard inequality on harmonic maps, that will still work in our situation, we will obtain a Schwarz lemma controlling uniformly the gradient of the disks. And we conclude using the following lemma used already in \cite{iv04} and \cite{hg25}.

\begin{Lem}\label{lem:1}
	\hspace{\parindent} Let $\Mc$ a Kobayashi hyperbolic Riemannian manifold, $\Sigma$ a hypersurface and $p\in\Sigma$. If there exists a $\Cc^1$ map $\delta$ from $\Mc$ to $\RR$ such that $\delta(p)=0$, $\delta_{\vert\Mc-\Sigma}\neq0$ and $C>0$ such that for all minimal disks $u:\DD\rightarrow\Mc-\Sigma,\vert\nabla(\delta\circ u)\vert(0)\leq C\vert\delta\circ u\vert(0)$ then $p$ is at infinite distance from $\Mc-\Sigma$.
\end{Lem}

\begin{Not}
	Here $\vert\cdot\vert$ denotes the absolute value and the Euclidean norm. Later in this paper when we will require to use the metric we will write $\vert\cdot\vert_g$ for the $g$-norm of a tangent vector.
\end{Not}

\begin{proof}[Proof of Lemma \ref{lem:1}]
	Let $q\in\Mc-\Sigma$, $p'$ that needs to be thought as near $p$ and $\gamma$ a $\Cc^1$ path from $p'$ to $q$. Note that $\Sigma$ being a hypersurface, we can suppose that $\delta\circ\gamma> 0$ without loss of generality. By definition, the length of $\gamma$ is $\ell_{\Mc-\Sigma}(\gamma):=\displaystyle\int_0^1 k_{\Mc-\Sigma}(\gamma'(t))dt$. However, at a fixed $t\in[0,1]$, there exists a minimal disk $u^t:(x,y)\in\DD\mapsto u^t(x,y)\in\Mc-\Sigma$ with $u^t(0)=\gamma(t)$ and $d(u^t)_0(1)=\dfrac{1}{2k_{\Mc-\Sigma}(\gamma'(t))}\gamma'(t)$. From this we get:
	\begin{align*}
	\dfrac{\partial(\delta\circ u^t)}{\partial x}(0) &= d\delta_{\gamma(t)}\left(\dfrac{\partial u^t}{\partial x}(0)\right)\\[1em]
	&= \dfrac{1}{2k_{\Mc-\Sigma}(\gamma'(t))}d\delta_{\gamma(t)}(\gamma'(t))\\[1em]
	&= \dfrac{1}{2k_{\Mc-\Sigma}(\gamma'(t))}(\delta\circ\gamma)'(t).
	\end{align*}
	Therefore, because $k_{\Mc-\Sigma}(\gamma'(t))>0$:
	\begin{align*}
	k_{\Mc-\Sigma}(\gamma'(t)) &= \dfrac{1}{2}\left\vert\dfrac{(\delta\circ\gamma)'(t)}{\partial_x(\delta\circ u^t)(0)}\right\vert \\[1em]
	&\geq \dfrac{1}{2C}\dfrac{\vert(\delta\circ\gamma)'(t)\vert}{\vert\delta\circ u^t(0)\vert} \mbox{ thanks to the control hypothesis on gradients,}\\[1em]
	&= \dfrac{1}{2C}\dfrac{\vert(\delta\circ\gamma)'(t)\vert}{\delta\circ\gamma(t)} \\[1em]
	&\geq \dfrac{1}{2C}\dfrac{(\delta\circ\gamma)'(t)}{\delta\circ\gamma(t)}.
	\end{align*}
	In particular:
	\begin{align*}
	\ell_{\Mc-\Sigma}(\gamma) &= \displaystyle\int_0^1k_{\Mc-\Sigma}(\gamma'(t))dt \\[1em]
	&\geq \dfrac{1}{2C}\displaystyle\int_0^1\dfrac{(\delta\circ\gamma)'(t)}{\delta\circ\gamma(t)}dt \\[1em]
	&= \dfrac{1}{2C}\displaystyle\int_{\delta(p')}^{\delta(q)}\dfrac{ds}{s} \\[1em]
	&= \dfrac{1}{2C}(\log\delta(q)-\log\delta(p'))\xlongrightarrow[p'\rightarrow p]{}=+\infty.
	\end{align*}
	So $p$ is at infinite distance.
\end{proof}

\subsection{Schwarz lemma on uniformly \texorpdfstring{$\Delta$}{Laplacian}-controlled map}

\hspace{\parindent} Let's first establish the Schwarz lemma. The trick is going to use a Harnack inequality. It states usually that, given a harmonic map $u:\DD\rightarrow\RR_+$, we have $\sup u\leq C\inf u$ where $C$ is independent of $u$ and the infimum and supremum are taken over a relatively compact domain in the disk. A similar result stay true in a broader case where $\vert\Delta u\vert\leq Cu$, due to results from \cite{tr67}.

Note that part of this subsection could be generalized to quasi-linear elliptic operators because the weak Harnack inequalities used in the proof stays true. However the Schwarz lemma requires an inequality on the Poisson equation which needs the standard Laplacian.

\par Let $\Omega$ be a domain of $\RR^n$ that we can suppose to be a neighborhood of $0$ without loss of generality. Let $\Cc^\infty(\Omega)$ be the set of smooth maps from $\Omega$ to $\RR$ and $\Cc^\infty_+(\Omega)$ be the subset of smooth non-negative maps. Let $Q(\rho)$ be the ball in $\Omega$ around $0$ of radius $\rho$ for the norm $\vert v\vert_\infty:=\max\vert v^i\vert$ in $\RR^n$. And finally we write $\vert\cdot\vert_{L^\gamma(\Omega)}$ for the usual $L^\gamma$-norm,$\gamma>1$.

\begin{Def}[$\Delta$-controlled maps]
	\hspace{\parindent} A set $\mathfrak{U}\subset\Cc^\infty(\Omega)$ is said to be $\Delta$-controlled if there exists a constant $C_\mathfrak{U}\geq0$ such that for every $u\in\mathfrak{U},\vert\Delta u\vert\leq C_\mathfrak{U}\vert u\vert$.
\end{Def}

\begin{Rmk}
	In particular, $\Delta$-controlled maps with small absolute value can be considered as \textit{almost} harmonic maps. It will be the case for a transversal component of minimal disks approaching $\Sigma$ in the proof of Theorem \ref{thm:1}. Moreover, one can hope that for \textit{almost} harmonic maps some results stay true, which is the case for the Harnack inequality.
\end{Rmk}

We will look at subsolutions and supersolutions of the differential equation:
\begin{equation}\label{eq:1}
(\Delta+C)u=0
\end{equation}
where $C$ is a constant. We say that $u$ is a supersolution of Equation \eqref{eq:1} if $(\Delta+C)u\leq0$ and a subsolution if $(\Delta+C)u\geq0$. Note that a map in a $\Delta$-controlled set is both a supersolution of $(\Delta-C_\mathfrak{U})u=0$ and a subsolution of $(\Delta+C_\mathfrak{U})u=0$.

\begin{Pro}[Harnack inequality]\label{pro:h}
	\hspace{\parindent} Let $\mathfrak{U}$ be a $\Delta$-controlled set and $\rho>0$ such that $Q(3\rho)$ is compactly contained in $\Omega$. There is a constant $C'$ such that:
	$$\forall u\in\mathfrak{U},\sup_{Q(\rho)}u\leq C'\inf_{Q(\rho)}u$$
	with $C'$ depending only on the dimension $n$ of the domain $\Omega$, the constant $C_\mathfrak{U}$ and the radius $\rho$.
\end{Pro}

To start, note that Trudinger already gave a Harnack inequality for solutions of Equation \eqref{eq:1} in \cite{tr67}. The proof is rather immediate once he established the two weak Harnack inequalities given below. See \cite{tr67} for the proof of both inequalities.

\begin{Lem}[Weak Harnack inequality for supersolutions, Trudinger \cite{tr67}]\label{lem:wk1}
	\hspace{\parindent} Let $u$ be a non-negative bounded supersolution of Equation \eqref{eq:1} on $Q(3\rho)$. Then there is a constant $C_+\geq 0$ such that:
	$$\dfrac{1}{\rho^{n/\gamma}}\vert u\vert_{L^\gamma(Q(2\rho))}\leq C_+\inf_{Q(\rho)}u$$
	for any $\gamma<n/(n-2)$ if $n\neq2$ or any $\gamma>1$ when $n=2$, with $C_+$ depending on $n$, $C$, $\rho$ and $\gamma$.
\end{Lem}

\begin{Lem}[Weak Harnack inequality for subsolutions, Trudinger \cite{tr67}]\label{lem:wk2}
	\hspace{\parindent} Let $u$ be a non-negative bounded subsolution of Equation \eqref{eq:1} on $Q(3\rho)$. Then there is a constant $C_-\geq 0$ such that:
	$$\sup_{Q(\rho)}u\leq\dfrac{C_-}{\rho^{n/\gamma}}\vert u\vert_{L^\gamma(Q(2\rho))}$$
	for any $\gamma>1$, with $C_+$ depending on $n$, $C$, $\rho$ and $\gamma$.
\end{Lem}

Now, to prove the Harnack inequality for solutions of Equation \eqref{eq:1}, one just has to combine both weak Harnack inequalities. This is possible because a solution of Equation \eqref{eq:1} is both a subsolution and a supersolution of the equation. Therefore with both Lemma \ref{lem:wk1} and Lemma \ref{lem:wk2}, one obtains:
$$\sup_{Q(\rho)}u\leq\dfrac{C_-}{\rho^{n/\gamma}}\vert u\vert_{L^\gamma(Q(2\rho))}\leq C_-C_+\inf_{Q(\rho)}u,$$
by choosing $\gamma$ between $1$ and $n/(n-2)$ ($=\infty$ when $n=2$).

The difference in our situtation, is that a $\Delta$-controlled map is not a solution of Equation \eqref{eq:1}. However, the proof persists because a $\Delta$-controlled map is still both subsolution and supersolution of two different equations of the same type of Equation \eqref{eq:1}. In particular, we can still apply both lemmas.

\begin{proof}[Proof of Proposition \ref{pro:h}]
Firstly, a $\Delta$-controlled map $u$ is a supersolution of $(\Delta-C_\mathfrak{U})u=0$. Therefore, by Lemma \ref{lem:wk1} we have
$$\dfrac{1}{\rho^{n/\gamma}}\vert u\vert_{L^\gamma(Q(2\rho))}\leq C_+\inf_{Q(\rho)}u,\forall\gamma<n/(n-2).$$
Secondly, $u$ is subsolution of $(\Delta+C_\mathfrak{U})u=0$. Lemma \ref{lem:wk2} gives us:
$$\sup_{Q(\rho)}u\leq\dfrac{C_-}{\rho^{n/\gamma}}\vert u\vert_{L^\gamma(Q(2\rho))},\forall\gamma>1.$$
Finaly, fixing $1<\gamma<n/(n-2)$ we get:
$$\sup_{Q(\rho)}u\leq\dfrac{C_-}{\rho^{n/\gamma}}\vert u\vert_{L^\gamma(Q(2\rho))}\leq C_-C_+\inf_{Q(\rho)}u.$$
\end{proof}

The important point is that the constant appearing in the inequality does not depend on the choice of map $u$. In particular we can derive the following Schwarz lemma.

\begin{Lem}[Schwarz lemma for $\Delta$-controlled map]\label{lem:2}
	\hspace{\parindent} Let $\mathfrak{U}\subset\Cc^\infty_+(\Omega)$ be a $\Delta$-controlled set. Let $\rho$ such that $Q(3\rho)$ is compactly contained in $\Omega$. Then there is a constant $C''$ depending only on $n$,$\rho$ and $C_\mathfrak{U}$ such that:
	$$\forall u\in\mathfrak{U},\vert\nabla u\vert(0)\leq C''u(0).$$
\end{Lem}

\begin{proof}
	To start, there is a well known gradient estimate for Poisson problem solutions. See \cite{tr83}. Given $u:\overline{Q(\rho)}\rightarrow\RR$ solution of the Poisson problem $\Delta u=f$ with $u\in\Cc^2(Q(\rho))\cap\Cc^0(\overline{Q(\rho)})$ and $f$ bounded in $Q(\rho)$ then:
	$$\forall i\in[\![1,n]\!],\left\vert\dfrac{\partial u}{\partial x^i}\right\vert(0)\leq\dfrac{n}{\rho}\sup_{\partial Q(\rho)}\vert u\vert+\dfrac{\rho}{2}\sup_{Q(\rho)}\vert f\vert.$$
	In particular, having $\sup_{\partial Q(\rho)}\vert u\vert\leq\sup_{Q(\rho)}\vert u\vert$ by continuity of $u$ and summing over all the indices we get that:
	$$\vert\nabla u\vert(0)\leq\dfrac{n^2}{\rho}\sup_{Q(\rho)}\vert u\vert+\dfrac{\rho n}{2}\sup_{Q(\rho)}\vert f\vert.$$
	Now taking $u\in\mathfrak{U}$, $u$ is solution of the trivial Poisson problem where $f:=\Delta u$, therefore:
	$$\vert\nabla u\vert(0)\leq\dfrac{n^2}{\rho}\sup_{Q(\rho)}u+\dfrac{\rho n}{2}\sup_{Q(\rho)}\vert\Delta u\vert.$$
	Moreover by hypothesis:
	$$\vert\nabla u\vert(0)\leq\left(\dfrac{n^2}{\rho}+\dfrac{\rho nC_\mathfrak{U}}{2}\right)\sup_{Q(\rho)}u.$$
	Thus, using the Harnack inequality from Proposition \ref{pro:h} we have:
	$$\vert\nabla u\vert(0)\leq\left(\dfrac{n^2}{\rho}+\dfrac{\rho nC_\mathfrak{U}}{2}\right)\sup_{Q(\rho)}u\leq\left(\dfrac{n^2}{\rho}+\dfrac{\rho nC_\mathfrak{U}}{2}\right)C'\inf_{Q(\rho)}u\leq\left(\dfrac{n^2}{\rho}+\dfrac{\rho nC_\mathfrak{U}}{2}\right)C'u(0).$$
\end{proof}

Now, for the remaining part of the proof, we only need to find an appropriate coordinates system such that a transversal component of minimal disks will be $\Delta$-controlled. With such control, one will be able to obtain Theorem \ref{thm:1} by combining both Lemma \ref{lem:1} and Lemma \ref{lem:2}.

\subsection{Distance to totally geodesic hypersurfaces}

\hspace{\parindent} As stated above, the first goal is to find a suitable coordinates system such that one coordinate component of minimal disks are \textit{almost harmonic} in the above sense of $\Delta$-controlled maps. Then we will be able to use the Schwarz lemma for the gradient control. Let's start by constructing the coordinates.

Let $N\Sigma:=\{(p,v),q\in\Sigma,v\in T_p\Sigma^\bot\}$ be the normal bundle over $\Sigma$. One can show that the exponential map (restricted to the normal bundle) is a diffeomorphism onto its image locally around the zero section. In particular, we define $\Oc_\Sigma$ the largest open neighborhood of $0$ in $N\Sigma$ such that $\exp_{\vert\Oc_\Sigma}$ is a diffeomorphism onto its image.

Now given a coordinate system $(y^1,...,y^{n-1})$ of $\Sigma$ in a neighborhood $\Vc$ of a point $p\in\Sigma$ and an orthonormal section $\nu$ of $N\Sigma$ we define:
$$\forall i\in[\![1,n-1]\!],\forall q\in\Vc,x^i\left(\exp_{\vert\Oc_\Sigma}\left(s\nu(q)\right)\right)=y^i(q)$$
$$\forall q\in\Vc,t\left(\exp_{\vert\Oc_\Sigma}\left(s\nu(q)\right)\right)=s$$
provided that $s$ is small enough for $s\nu\in\Oc_\Sigma$.

Since $\exp_{\vert\Oc_\Sigma}$ is a diffeomorphism, the above equations defines, in a neighborhood of $p$, a coordinate system $(x^1,...,x^{n-1},t)$. This is actually a special case of Fermi's coordinates of a submanifold where the submanifold is a hypersurface. See \cite{gr03}.

\begin{Rmk}
	A few important properties:
\begin{itemize}
	\item[--] The normal coordinate $t$ is the signed distance to $\Sigma$. This is easily shown by the definition above as $t$ is obtained using normal geodesics.
	\item[--] The gradient $\nabla t$ of the normal coordinate is the unit normal field of each leaf in the foliation $(\Sigma_t)_t$ of level sets of $t$. Indeed if $X$ is tangential to the level set $\Sigma_c:=t^{-1}(\{c\})$ then $dt(X)=0$ because $T_q\Sigma_c=\ker dt_q,q\in\Sigma_c$. It tells us that $\nabla t(q)$ is normal to $\Sigma_{t(q)}$. Moreover, given $\gamma(s)=\exp_q(s\nu(q))$ a normal geodesic, where $q\in\Sigma\cap\Uc$, then by construction $t(\gamma(s))=s$. Therefore:
\begin{align*}
1 &\equiv \dfrac{dt(\gamma(s))}{ds} \\[1em]
&= dt_{\gamma(s)}(\gamma'(s)) \\[1em]
&= dt_{\gamma(s)}(\nu(\gamma(s))) \\[1em]
&= g_{\gamma(s)}(\nabla t(\gamma(s)),\nu(\gamma(s))).
\end{align*}
Since this is true for all normal geodesic $\gamma$, then $1\equiv g(\nabla t,\nu)$. Having $\nabla t$ normal to the level set, we get that $\nabla t$ is the normal vector field of the level sets. In particular, having $\nabla t_{\vert\Sigma}=\nu$, we will still write $\nu$ for $\nabla t$.
	\item[--] The extended normal field $\nu$ is a geodesic vector field that is $\nabla_\nu\nu=0$. Indeed using the same notations as above $\nabla_\nu\nu=\nabla_{\gamma'}\gamma'=0$ because $\gamma$ is a geodesic.
\end{itemize}
\end{Rmk}

Now that we have an adapted coordinates system, we want to show that, given a disk $u:\DD\rightarrow\Mc-\Sigma$, $t\circ u$ is uniformly $\Delta$-controlled to be able to use Lemma \ref{lem:2}. To do that we first need to compute the Hessian of $t$, $\Hess t(X,Y):=g(\nabla_X\nabla t,Y)=g(\nabla_X\nu,Y),X,Y\in\Gamma(T\Uc)$. Immediately, if $X$ and $Y$ are tangential to the level sets, meaning that for every $q\in\Uc,X(q),Y(q)\in T_q\Sigma_{t(q)}$ then $\Hess t(X,Y)=-\mbox{I}\!\mbox{I}_{\Sigma_t}(X,Y)$ where $\mbox{I}\!\mbox{I}_{\Sigma_t}$ is the second fundamental form of the level set of $\Sigma_t$. Moreover, having $\nabla_\nu\nu=0$, the tangential/normal and normal/normal components of the Hessian are trivial. In particular, the Hessian can be written as $\Hess t=\displaystyle\left(\begin{array}{cc}
-\mbox{I}\!\mbox{I}_{\Sigma_t} & 0 \\ 
0 & 0
\end{array} \right)$.

Let's now control the way the second fundamental forms of the level sets evolve approaching $\Sigma$. To do this we will first control the shape operators.

\begin{Not}
	To control the shape operators, we are going to use the Frobenius norm that we will write $\vert\cdot\vert_g$ as it extends the Riemannian norm of vectors that we already wrote that way. Given a $(p,q)$-tensor $T$, the Frobenius norm of $T$ is $\vert T\vert_g:=\sqrt{g_{i_1i'_1}...g_{i_pi'_p}g^{j_1j'_1}...g^{j_qj'_q}T_{j_1...j_q}^{i_1...i_p}T_{j'_1...j'_q}^{i'_1...i'_p}}$.
\end{Not}
	
	Finally, the last notion we need is parallel transport. Given a path $\gamma:I\rightarrow\Mc$, it defines a set of smooth isometries $\{\Gamma_s^t:T_p^q\Mc_{\gamma(t)}\rightarrow T_p^q\Mc_{\gamma(s)},s,t\in I\}$ called parallel transport on $(p,q)$-tensor. In particular, the parallel transport can recover the connection with the following formula:
$$\forall T\in\Gamma(T_p^q\Mc),(\nabla_{\gamma'}T)_{s}=\left.\dfrac{d}{dt}\Gamma_s^t T_{\gamma(t)}\right\vert_{t=s}.$$
For more details on these notions, see \cite{lee6}. Note that we are going to use an equivalent formula. Let us fix $t_0\in I$ and $\tilde{T}(t)=\Gamma_{t_0}^t T_{\gamma(t)}$. From standard properties on the parallel transport, we have that $\tilde{T}(t)=\Gamma_{t_0}^s\Gamma_s^t T_{\gamma(t)}$. Now by using the previous formula we get:
\begin{align*}
	\left.\dfrac{d}{dt}\tilde{T}(t)\right\vert_{t=s}&=\Gamma_{t_0}^s\left(\left.\dfrac{d}{dt}\Gamma_s^t T_{\gamma(t)}\right\vert_{t=s}\right) \\[1em]
	&= \Gamma_{t_0}^s(\nabla_{\gamma'}T)_{s}.
\end{align*}
This being true for all $s\in I$, it gives us $\vert \tilde{T}'\vert_g=\vert\nabla_{\gamma'}T\vert_g$ in Frobenius norm because the parallel transport maps are isometries. We will use this last equality in the proof of the Theorem \ref{thm:1}.

\begin{Lem}\label{lem:abv}
	\hspace{\parindent} Let $\Mc$ be a Riemannian manifold, $\Sigma$ a hypersurface and $\Uc$ a sufficiently small neighborhood of a point of $\Sigma$. Then $\Sigma\cap\Uc$ is totally geodesic if and only if it exists a real constant $C$ such that for every $q\in\Uc,\vert S_{q}\vert_{g'}\leq Cd(q,\Sigma)$ where $S_q$ is the shape operator at $q$ of the unique leaf containing $q$ and $g'$ is the first fundamental form of the leaf.
\end{Lem}

\begin{Rmk}
	To be precise, by $\Uc$ small enough we mean that $\Uc$ must be compactly contained in the domain of Fermi's coordinates.
\end{Rmk}

\begin{proof}
	Firstly, the $\Leftarrow$ implication is easy. Indeed assume the control inequality on the shape operators is true. Let $q\in\Sigma\cap\Uc$ and $\gamma$ be a normal geodesic starting from $q$. Note that $s\mapsto\vert S_{\gamma(s)}\vert_g$ is continuous. Indeed the shape operators can be viewed as a $(1,1)$-tensor $S$ which encodes all data of the shape operators. To be precise, let $E$ be the tangent bundle of the foliation $(\Sigma_t)_t$, then $S$ can be seen as a section of $T_1^1E=E\otimes E^*$. In particular $x\mapsto S_x$ is continuous, so is $x\mapsto\vert S_x\vert_{g'}$ and $s\mapsto\vert S_{\gamma(s)}\vert_{g'}$. In particular, by taking the limit in the inequality of Lemma \ref{lem:abv}, it shows that $\vert S_q\vert_{g'}\leq0$. This being true for all base point $q\in\Sigma\cap\Uc$, we conclude that $\Sigma\cap\Uc$ is totally geodesic.
	
	Secondly for the direct implication, it will be useful to see the shape operators as a global $(1,1)$-tensor on $\Uc$. Indeed, we can defined $S$ as a section of $T_1^1\Uc$ by looking at $X\mapsto\nabla_X\nu$ where $\nu$ is the normal vector field to the foliation. Note here that $S$ can be restricted on the set of tangent vector fields to a given leaf, therefore $\vert S_{\vert\Sigma_t}\vert_{g'}$ is well defined given a leaf $\Sigma_t$. Moreover, by definition of the Frobenius norm $\vert S_{\vert\Sigma_t}\vert_{g'}\leq\vert S\vert_g$ on $\Sigma_t\cap\Uc$. Therefore to show the intended inequality, we just need to show a similar one to $\vert S\vert_g$.
	
	Note that since $\nu=\nabla t$, $\nu$ is a smooth vector field. Therefore $\nabla_\nu S$ is a smooth $(1,1)$-tensor, in particular $\vert\nabla_\nu S\vert_g$ is continuous where it is defined that is on the domain of Fermi's coordinates. Because $\Uc$ is compactly contained in the domain of Fermi's coordinates, there is constant $C$ such that $\vert\nabla_\nu S\vert_g\leq C$ on $\Uc$. Moreover given $\gamma$ a normal geodesic, it defines a parallel transport $\{\Gamma_s^t\}$ as stated in the previous remark. Let $\tilde{S}(s):=\Gamma^s_0(S_{\gamma(s)})$. It gives us:
	$$\vert\tilde{S}'\vert_g=\vert\nabla_{\gamma'}S\vert_g=\vert\nabla_\nu S\vert_g\leq C.$$
	Hence, $\forall s,\vert \tilde{S}(s)\vert_g\leq Cs$, because $\tilde{S}(0)\equiv 0$. Moreover the parallel transport is an isometry. Therefore $\vert S_{\gamma(s)}\vert_g=\vert\tilde{S}(s)\vert_g\leq Cs$. This being true for all normal geodesics $\gamma$, it gives us the intended inequality.	
\end{proof}

In particular, we get the following corollary.

\begin{Cor}\label{cor:1}
	\hspace{\parindent} Let $\Mc$ be a Riemannian manifold, $\Sigma$ a hypersurface and $\Uc$ a small enough neighborhood of a point of $\Sigma$. Then  $\Sigma\cap\Uc$ is totally geodesic if and only if it exists $C\in\RR$ such that for every $q\in\Uc,\vert(\mbox{\emph{I}}\!\mbox{\emph{I}}_{\Sigma_{t(q)}})_q\vert_{g'}\leq Cd(q,\Sigma)$.
\end{Cor}

We have now every tool to prove Theorem \ref{thm:1}.

\begin{Thm*}
	\hspace{\parindent}Let $(\Mc,g)$ be a Kobayashi hyperbolic Riemannian manifold, $\Sigma$ a hypersurface and $p\in\Sigma$. If $\Sigma$ is locally totally geodesic around $p$ then $p$ is at infinite Kobayashi distance.	
\end{Thm*}

\begin{proof}
	Let $\Uc$ be a relatively compact open neighborhood of $p$ in $\Mc$ small enough to assure the existence of Fermi's coordinates, for $\Uc\cap\Sigma$ to be totally geodesic and small enough to be able to apply Corollary \ref{cor:1}. To keep similar notations $t:\Uc\rightarrow\RR$ will be the normal coordinate. Note that $\Sigma$ cuts $\Uc$ into two connected components $\{t>0\}$ and $\{t<0\}$. Therefore, when taking a disk $u:\DD\rightarrow\Uc-\Sigma$, we can always suppose that $t\circ u>0$ without loss of generality. It suffices to show that the set $\{u^\bot:=t\circ u>0,u:\DD\xrightarrow{\mbox{harm.}}\Uc\}$ is uniformly $\Delta$-controlled. Indeed, with this hypothesis, we can apply the Schwarz lemma in Lemma \ref{lem:2} to finally conclude with Lemma \ref{lem:1}.
	
	Let $u:\DD\rightarrow\Uc$ be a harmonic disk. Firstly,
	\begin{align*}
	\Hess u^\bot(X,Y)&=\Hess t\circ u(X,Y) \\[1em]
	&= X\cdot(Y\cdot(t\circ u))-\nabla_XY\cdot(t\circ u) \\[1em]
	&= X\cdot(dt(du(Y))-dt(du(\nabla_XY)) \\[1em]
	&= \left(\nabla_{du(X)}dt\right)(du(Y))+dt\left(\nabla_{du(X)}du(Y)\right)-dt\left(du(\nabla_XY)\right) \\[1em]
	&= \Hess t(du(X),du(Y))+dt\left(\nabla_{du(X)}du(Y)-du(\nabla_XY)\right) \\[1em]
	&= \left(u^*\Hess t\right)(X,Y)+dt((\nabla du)(X,Y)).
	\end{align*}
	Therefore $\Delta u^\bot:=\tr\Hess u^\bot=\tr\left(u^*\Hess t\right)+dt(\tau(u))$ where $\tau(u)$ is the tension field of $u$. But $u$ is harmonic, so $\tau(u)=0$, hence $\Delta u^\bot=\tr u^*\Hess t$.
	
	Secondly, we claim that $\vert\Delta u^\bot\vert\leq\sqrt{2}\Vert du\Vert\vert\Hess t\vert_g$ where $\Vert\cdot\Vert$ is the operator norm and $\vert\cdot\vert_g$ is the Frobenius norm on tensors relative to the metric $g$. To show that, let $T$ be a $(0,2)$-tensor on the disk $\DD$. The Frobenius norm of $T$ for the Euclidean metric $g_0$ using the canonical coordinates is $\vert T\vert_{g_0}=\sqrt{\sum_{i,j}T_{ij}^2}$. Moreover, we have:
	$$\begin{array}{rl}
	\vert\tr T\vert&=\vert\sum_i T_{ii}\vert \\[1em]
	&\leq \sqrt{2}\sqrt{\sum_{i}T_{ii}^2}\mbox{ using Cauchy-Schwarz inequality,} \\[1em]
	&\leq \sqrt{2}\sqrt{\sum_{i,j}T_{ij}^2} \\[1em]
	&\leq \sqrt{2}\vert T\vert_{g_0}.
	\end{array}$$
	In particular, at a point $z\in\DD$, we have:
	\begin{align*}
	\vert\Delta u^\bot(z)\vert&\leq\sqrt{2}\vert u^*\Hess t\vert_{g_0}(z) \\[1em]
	&\leq\sqrt{2}\Vert du_z\Vert\vert\Hess t_{u(z)}\vert_g
	\end{align*}
	which is the intended inequality.
	
	Now, recall that $\Hess t=\displaystyle\left(\begin{array}{cc}
-\mbox{I}\!\mbox{I}_{\Sigma_t} & 0 \\ 
0 & 0
\end{array} \right)$. Therefore, Corollary \ref{cor:1} gives us a control on the Hessian $\vert\Hess t\vert_g\leq Ct$. In particular, it shows us that $\vert\Delta u^\bot\vert\leq\sqrt{2}C\Vert du\Vert t\circ u$. Note that $C$ is independent of $u$, this is extremely important to apply Lemma \ref{lem:2}.

	Finally, we have to control uniformly $\Vert du\Vert$. Up to reducing $\Uc$, by Royden hyperbolicity of $\Mc$, there is a constant $c>0$ such that for all $(q,v)\in T\Uc,k_\Mc(q,v)\geq c\vert v\vert_g$. In particular, applying this inequality for $(u(z),du_z(\xi))$, with $(z,\xi)\in\DD\times\CC$ gives us:
	\begin{equation}\label{eq:2}
		\vert du_z(\xi)\vert_g\leq\dfrac{1}{c}k_{\Mc}(u(z),du_z(\xi))\leq\dfrac{1}{c}\rho_{\DD}(z,\xi)\leq\dfrac{1}{c}\dfrac{\vert\xi\vert}{1-\vert z\vert^2}.
	\end{equation}
Now, let us fix a radius $\rho$ small enough for $Q(3\rho)$ to be compactly contained in $\DD$. There are two reasons for us to do so. Firstly, restricting $z$ to $Q(3\rho)$ gives us a uniform control of the $1/(1-\vert z\vert^2)$ term in Equation \eqref{eq:2}. Secondly, it is one of the hypothesis to apply Lemma \ref{lem:2}. Therefore, we have a uniform constant $c'$ such that:
$$\forall(z,\xi)\in Q(3\rho)\times\CC,\vert du_z(\xi)\vert_g\leq c'\vert\xi\vert.$$ This last inequality gives us, by definition of the operator norm, that $\Vert du\Vert\leq c'$ on $Q(3\rho)$.

	Finally we obtained a uniform constant $C'$ such that $\vert\Delta u^\bot\vert\leq C' t\circ u$ on $Q(3\rho)$. Therefore, by Lemma \ref{lem:2}, there is a uniform constant $C_1$ such that $\vert\nabla u^\bot\vert(0)\leq C_1u^\bot(0)$ which can be written as $\vert\nabla\delta\circ u\vert(0)\leq C_1\delta\circ u(0)$ with $\delta$ the signed distance to $\Sigma$. We can conclude, with Lemma \ref{lem:1}, that $p$ is at infinite distance.
\end{proof}

\section{Proof of Theorem 2}
\label{sec:3}

\hspace{\parindent} The key idea to prove Theorem \ref{thm:2} is to find an explicit sequence of points going to $\Sigma$ with bounded distance. To do that we will first show two lemmas that will be important. Let us first reduce the situation. We claim that proving Theorem \ref{thm:2} on $\HH^n,\RR^n$ and $\SS^n$ is enough. Indeed, let $\kappa\in\RR$ and $\Mc_\kappa^n$ the associated model space, that is the unique complete simply connected Riemannian n-manifold of constant sectional curvature $\kappa$. As mentionned in the introduction, it is enough to show Theorem \ref{thm:2} on Kobayashi hyperbolic domains of model manifolds, because every manifold of constant sectional curvature is locally isometric to an open subset of $\Mc_\kappa^n$. Secondly, the Riemannian manifold $\Mc_\kappa^n$ is homothetic to $\HH^n$ if $\kappa<0$, $\RR^n$ if $\kappa=0$ of $\SS^n$ if $\kappa>0$ and homotheties preserve the geometry (in particular minimal surfaces). Hence solving the theorem for domains on these three spaces will give the general theorem.

In the proof, we will first study the situation in Kobayashi hyperbolic Euclidean domains then extend it on $\HH^n$ and $\SS^n$. The proofs are extremely similar (if not identical) so we will just add remarks on eventual changes.

\begin{Lem}\label{lem:ex}
\hspace{\parindent} Let $\Omega$ be a Kobayashi hyperbolic domain of $\Mc_\kappa^3$ a model $3$-manifold, $\Sigma\subset\Omega$ be a surface and $p\in\Sigma$. If $H_\Sigma(p)\neq 0$ where $H_\Sigma$ is the mean curvature of $\Sigma$, then the point $p$ is at finite distance from $\Omega-\Sigma$.
\end{Lem}

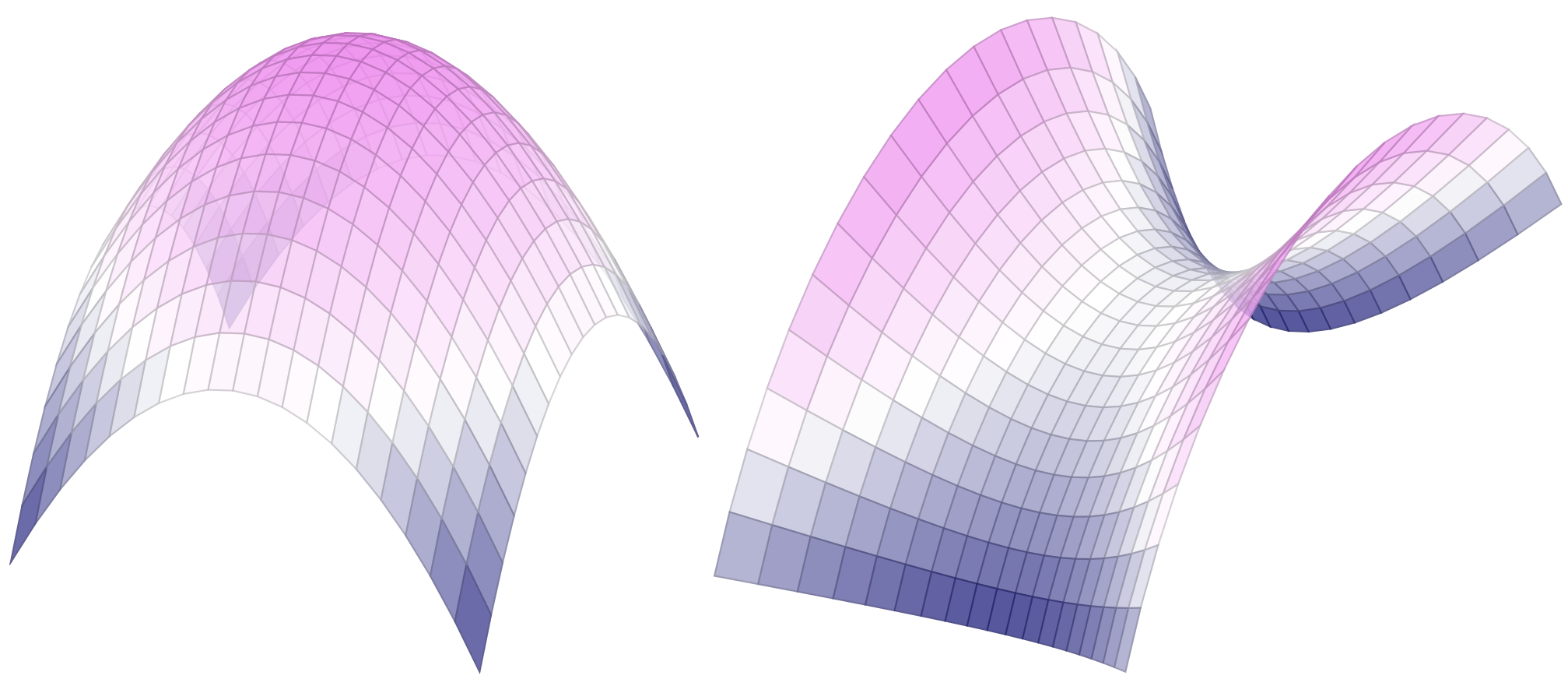
\begin{figure}[ht]
	\centering
	\begin{tikzpicture}
	\begin{axis}[axis equal, hide axis,scale=2]
		\addplot3[surf, colormap/violet, domain=-1:1, y domain=-1:1, samples=20, z buffer=sort, opacity=.8] ({x},{y},{-x*x-y*y});
	\end{axis}
	\end{tikzpicture}
	\begin{tikzpicture}
	\begin{axis}[axis equal, hide axis,scale=1.7]
		\addplot3[surf, colormap/violet, domain=-1:1, y domain=30:150, samples=20, z buffer=sort, opacity=.8] ({sinh(x)},{2 * cosh(x) * cos(y)},{2 * cosh(x) * sin(y)});
	\end{axis}
	\end{tikzpicture}
	\vspace{.5cm}
	\caption{The two types of surface.}
	\label{fig:sits}
\end{figure}

\begin{Rmk}
	The hypothesis $H_\Sigma(p)\neq 0$ corresponds to two situations:
	\begin{itemize}
		\item[--] either the two principal curvatures at $p$ have same sign and are non-zero, in particular $\Sigma$ is strictly convex in a neighborhood of $p$ or,
		\item[--] the two principal curvatures at $p$ have different signs and are different in absolute value.
	\end{itemize}
	See Figure \ref{fig:sits}. Note that by sign we mean positive, negative or zero. Therefore in the second situation, one principal curvature can vanish at $p$.
\end{Rmk}

\begin{proof}[Proof of Lemma \ref{lem:ex}]
	The first step of the proof is to find a strictly separated minimal surface to $\Sigma$ at the point of interest $p$. That is a minimal surface $S$ such that $p\in S$, $T_pS=T_p\Sigma$ and locally near $p$, $S-\{p\}\subseteq\Cc$ where $\Cc$ is one connected component of $\Omega-\Sigma$. In other words, $S$ stays strictly on one side of $\Sigma$ in a punctured neighborhood of $p$ and is tangent to $\Sigma$ at $p$.

	Let us first study the case of $\kappa=0$. Firstly, up to translations and rotations, we can suppose that $p=0$ and $T_0\Sigma=\{z=0\}$. Secondly, thanks to the implicit function theorem, we can find a parametrisation of $\Sigma$ in a neighborhood of $0$, of the form $(x,y)\mapsto(x,y,f(x,y))$. In particular, by the hypothesis on the tangent plane at $0$ we have that $f(x,y)=\dfrac{1}{2}(x\; y)\Hess f_{0}\;(x\;y)^T+o(x^2+y^2)$. Note that the Hessian being symmetric, thanks to the spectral theorem, up to a rotation in the plane $\{z=0\}$ we can suppose that the Hessian is diagonal. Therefore, without loss of generality, we can suppose that $\Sigma$ is the graph of a map $f:\DD\rightarrow\Omega$ such that
	\begin{itemize}
		\item[--] $f(x,y)=\pm\dfrac{1}{2}(k_xx^2+k_yy^2)+o(x^2+y^2),k_x,k_y>0$ or,
		\item[--] $f(x,y)=\dfrac{1}{2}(k_xx^2-k_yy^2)+o(x^2+y^2)$, $k_x\neq k_y,k_x\geq 0$ and $k_y>0$ (or $k_x>0$ and $k_y\geq 0$).
	\end{itemize}
	Note that $k_x$ and $k_y$ are the absolute values of the principal curvatures of $\Sigma$ at $0$.
	
	The surfaces we are going to use are catenoids. Recall that a catenoid is a surface given by the parametrisation $\Cc_k:(r,\theta)\mapsto\left(r,\dfrac{1}{k}\cosh(k r)\cos\theta,\dfrac{1}{k}\cosh(k r)\sin\theta\right),k>0$. This is a minimal surface and the principal curvatures at the neck $\subset\{x=0\}$ are $k$ and $-k$. After a translation, we can suppose that $0$ is on this circle and that $T_0\Cc_k=\{z=0\}$. In particular, locally near $0$ the catenoid can be seen as a graph of a map $g:\DD\rightarrow\Omega$ such that $g(x,y)=\dfrac{1}{2}(kx^2-ky^2)+o(x^2+y^2)$.
	
	The goal is to place a catenoid "above" (or "under" depending on the situation) $\Sigma$. With an adequate choice of $k$ the catenoid will be strictly separated. To do that, we just need to study the difference $g-f$.
	\begin{itemize}
		\item[--] Firstly, assume that the two principal curvatures have same sign. Let us work with negative curvatures (working with positive curvatures is similar). In particular, $(g-f)(x,y)=\dfrac{1}{2}((k_x+k)x^2+(k_y-k)y^2)+o(x^2+y^2)$. Taking $k<k_y$ gives us that $g-f>0$ in a punctured neighborhood of $0$.
		\item[--] Secondly, assume that the two principal curvatures have different signs and that $k_x\neq k_y$. Without loss of generality, we can assume that $k_x<k_y$. We have that $(g-f)(x,y)=\dfrac{1}{2}((k-k_x)x^2+(k_y-k)y^2)+o(x^2+y^2)$. Therefore, taking $k_x<k<k_y$ gives us that $g-f>0$ in a punctured neighborhood of $0$. See Figure \ref{fig:cat}.
	\end{itemize}
	Note that in the first situation, we could take $k=0$ which corresponds to the degenerate case of a plane.
	\begin{figure}[ht]
	\centering
	\begin{tikzpicture}
	\begin{axis}[axis equal, hide axis,scale=2]
		\addplot3[surf, colormap/violet, domain=-2:2, y domain=-2:2, samples=20, z buffer=sort, opacity=.8] ({y},{x},{(0.25 * x * x- y*y)/2});
		\addplot3[surf, colormap/hot, domain=0:1.4, y domain=0:360, samples=30, z buffer=sort, opacity=.8] ({x * sin(y)},{x * cos(y)},{(x * x)/4 * (cos(y) * cos(y) - sin(y) * sin(y))});
	\end{axis}
	\end{tikzpicture}
	\vspace{.5cm}
	\caption{A tangential piece of catenoid.}
	\label{fig:cat}
	\end{figure}
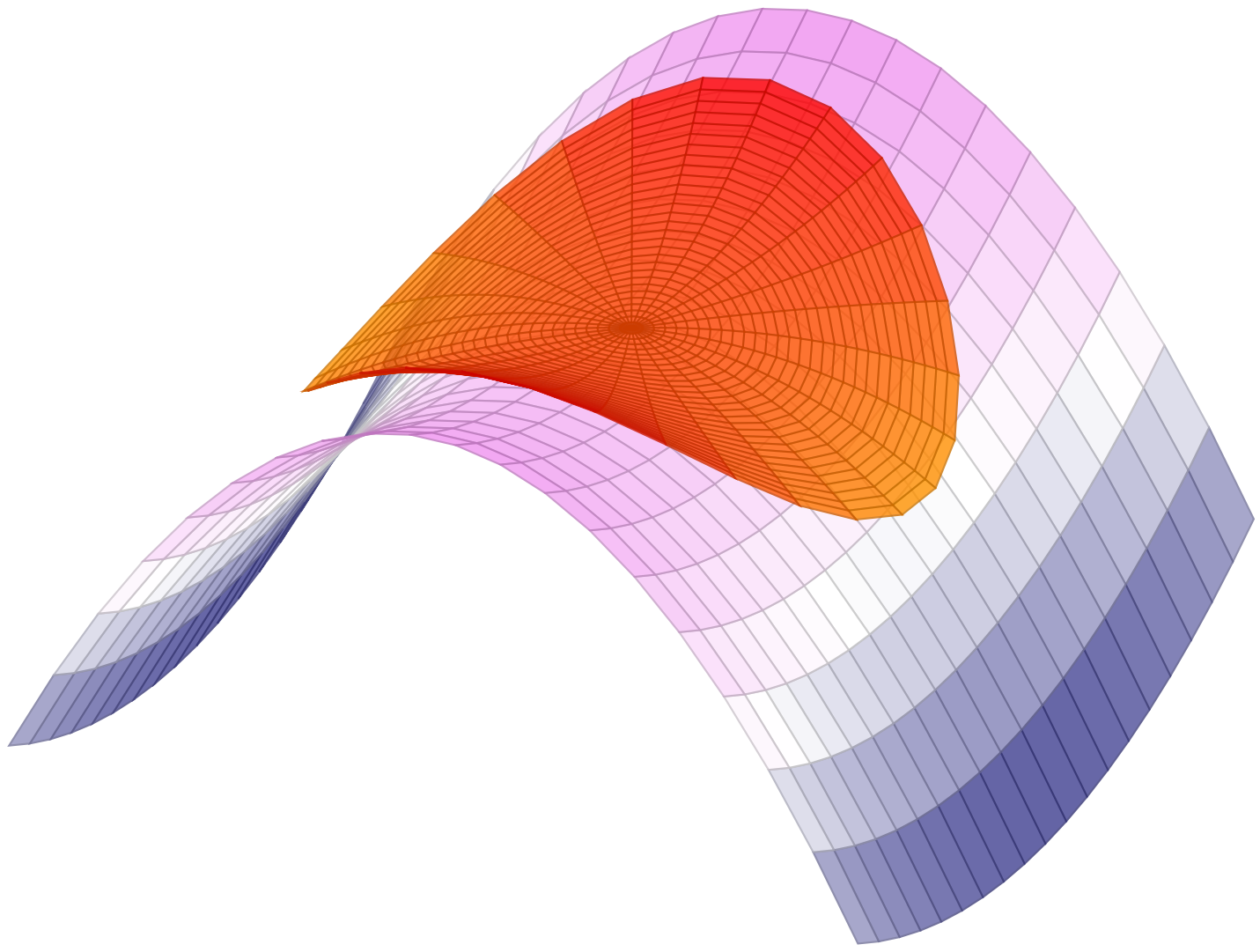	
	
	Let us now show this property on $\HH^n$ and $\SS^n$. We claim that it is enough to show that for all $k>0$ there exists a minimal surface in $\HH^n$ (respectively $\SS^n$) and a point $s$ of that minimal surface such that the principal curvature at $s$ are $\pm k$. Let us stress that on a model space $\Mc_\kappa^3$, there is a double transitivity from the isometry group. To be precise, the isometry group acts transitively
	\begin{itemize}
		\item[--] on the points of $\Mc_\kappa^3$, and
		\item[--] on the $2$-Grassmannians of $\Mc_\kappa^3$.
	\end{itemize}
	In particular, given a minimal surface with prescribed principal curvatures at a point, we can place it tangent to any $2$-plane in $T_p\Mc_\kappa^3$ for any $p$.
	
	Moreover, around the point $p\in\Sigma$, we can use normal coordinates and suppose that $p$ is sent to $0$ and $T_p\Sigma=\{z=0\}$ such that $\Sigma$ is a graph up to reducing the chart. By the previous paragraph, we can place a strictly separated minimal surface given the existence of one with prescribed principal curvatures because we are in the same situation as above.
	
	To show the existence of a minimal surface with prescribed principal curvatures at a point, we can use surfaces of revolution. See \cite{la70} and \cite{ca83}. Let us first do it for $\HH^3$. We define the $\alpha,\beta$-catenoid $\Cc_{\alpha,\beta}$, $\alpha,\beta>0$ by the immersion:
	$$f:(s,t)\mapsto(\cosh(s)\cosh(\alpha t),\cosh(s)\sinh(\alpha t),\sinh(s)\cos(\beta t),\sinh(s)\sin(\beta t)).$$
	Note here that $\HH^3$ is used with the hyperboloid model. A classical computation shows that the catenoid $\Cc_{\alpha,\beta}$ is a minimal surface and that the principal curvatures at the point $f(s,t)$ are $\pm\alpha\beta/\xi(s,t)$ where $\xi(s,t)=\alpha^2\cosh^2(s)+\beta^2\sinh^2(s)$. In particular, at the neck the principal curvatures are $\pm\beta/\alpha$, and $\beta/\alpha$ attains every value in $\RR_{>0}$ up to changing the parameters $\alpha$ and $\beta$. Now, as mentionned earlier, using the transitivity of the isometry group we can place this minimal surface "above" (or "under"). Then taking adequate parameters $\alpha$ and $\beta$, the minimal surface can be strictly separated in $\HH^3$ using the same curvature arguments as in the Euclidean case.
	
	The construction in $\SS^3$ is identical; with the immersion:
	$$f:(s,t)\mapsto(\cos(s)\cos(\alpha t),\cos(s)\sin(\alpha t),\sin(s)\cos(\beta t),\sin(s)\sin(\beta t)).$$
	The same computation shows that the image of $f$ is minimal and that the principal curvatures at $f(s,t)$ are $\pm\alpha\beta/\xi$ where $\xi(s,t)=\alpha^2\cos^2(s)+\beta^2\sin^2(s)$. In particular, in absolute value, the curvatures at the neck reach every value in $\RR_{>0}$ changing the parameters $\alpha$ and $\beta$. Again, by taking adequate parameters and using the transitivity we can construct a strictly separated minimal surface in $\SS^3$
	
	To summarize, we found a minimal disk $u_\infty:\DD\rightarrow\Omega$ such that $u_\infty(0)=p$, $u_\infty(\DD)$ is tangent to $\Sigma$ at $p$ and $u_\infty(\DD)$ stays strictly on one side of $\Sigma$ in a punctured neighborhood.
	
	Back to $\kappa=0$ to give the details. Let $p_k=(0,0,1/k)\in\Omega$ be a sequence of points going to $p$, and $u_k:\DD\rightarrow\Omega$ be the translated minimal disks such that $u_k(0)=p_k,\forall k$. See Figure \ref{fig:finite}. By construction, up to working locally, we can suppose that:
$$\forall k,u_k(\DD)\subset\Omega-\Sigma.$$
Moreover, because $u_\infty$ stays strictly (i.e. $u_\infty(\DD-\{0\})$) on one side of $\Sigma$, there exists a compact set of $\Omega-\Sigma$ containing $u_\infty(1/2)$. By construction, we can find a bigger compact set of $\Omega-\Sigma$ which contains all points $u_k(1/2)$ for every $k$. In particular, given a fixed point $q\in\Omega-\Sigma$ it gives us that $\sup_kd_{\Omega-\Sigma}(q,u_k(1/2))<\infty$.

	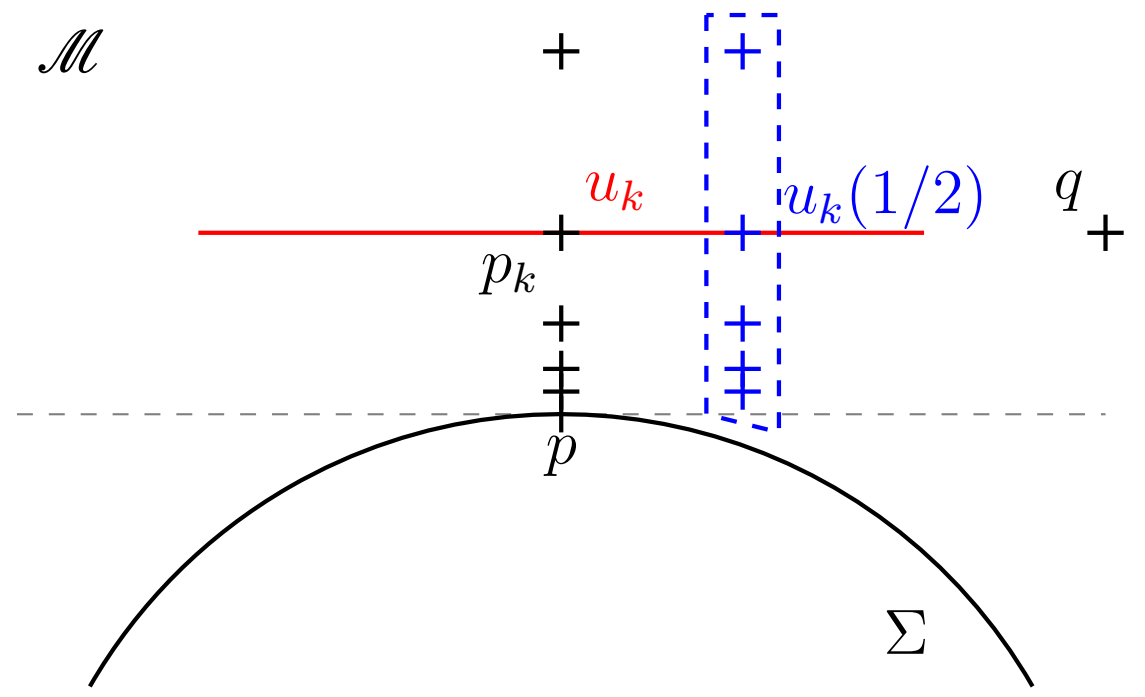
\begin{figure}[ht]
	\centering
	\begin{tikzpicture}[scale=1.2]
		\draw[dashed,gray] (-3,3) -- (3,3);
		\draw[thick] ([shift=(30:3)]0,0) arc (30:150:3);
		
		\draw[thick,red] (-2,4) -- (2,4) node[midway,anchor=south west] {$u_k$};
		
		\draw[thick] (0,2.9) -- (0,3.1);
		\draw (0,3) node[anchor=north] {$p$};
		
		\draw[thick] (0,4.9) -- (0,5.1);
		\draw[thick] (-.1,5) -- (.1,5);
		
		\draw[thick] (0,3.9) -- (0,4.1);
		\draw[thick] (-.1,4) -- (.1,4);
		\draw (0,4) node[anchor=north east] {$p_k$};
		
		\draw[thick] (0,3.4) -- (0,3.6);
		\draw[thick] (-.1,3.5) -- (.1,3.5);
		
		\draw[thick] (0,3.15) -- (0,3.35);
		\draw[thick] (-.1,3.25) -- (.1,3.25);
		
		\draw[thick] (0,3.025) -- (0,3.225);
		\draw[thick] (-.1,3.125) -- (.1,3.125);

		\draw[thick, blue] (1,4.9) -- (1,5.1);
		\draw[thick, blue] (.9,5) -- (1.1,5);
		
		\draw[thick, blue] (1,3.9) -- (1,4.1);
		\draw[thick, blue] (.9,4) -- (1.1,4);
		\draw[blue] (1.1,3.9) node[anchor=south west] {$u_k(1/2)$};
		
		\draw[thick, blue] (1,3.4) -- (1,3.6);
		\draw[thick, blue] (.9,3.5) -- (1.1,3.5);
		
		\draw[thick, blue] (1,3.15) -- (1,3.35);
		\draw[thick, blue] (.9,3.25) -- (1.1,3.25);
		
		\draw[thick, blue] (1,3.025) -- (1,3.225);
		\draw[thick, blue] (.9,3.125) -- (1.1,3.125);
		
		\draw[thick,dashed,blue] (0.8,5.2) -- (0.8,3) -- (1.2,2.9) -- (1.2,5.2) -- (0.8,5.2);

		\draw[thick] (3,3.9) -- (3,4.1);
		\draw[thick] (2.9,4) -- (3.1,4);
		\draw (3,4) node[anchor=south east] {$q$};
		
		\draw (1.9,1.8) node {$\Sigma$};
		\draw (-2.7,5) node {$\Mc$};
	\end{tikzpicture}
	\caption{A sequence of planar disks following the proof construction.}
	\label{fig:finite}
	\end{figure}
	
	Finally, $\forall k,d_{\Omega-\Sigma}(q,u_k(0))\leq\sup_kd_{\Omega-\Sigma}(q,u_k(1/2))+\rho_\DD(1/2,0)<\infty$. Therefore, $p$ is at finite distance from $\Omega-\Sigma$.
	
	To create such a sequence on domains in $\HH^3$ and $\SS^3$, we can replace the translations in the normal direction of $\Sigma$ at $p$ by loxodromic isometries on $\HH^3$ in the normal direction or rotations on $\SS^3$, again, in the normal direction. Both types of isometries will create (locally) adequate sequences. And the proof ends with the exact same arguments.
\end{proof}

\begin{Rmk}
	Note that, similar to the result of Ivashkovich-Rosay in \cite{iv04}, we cannot be more precise and look at the two connected components of $\Omega-\Sigma$. In particular, we can only say that the surface is at finite distance from one connected component. In some situation, the distance to a point can be finite on one component and infinite on the other. In particular when both principal curvatures have same sign. See Remark \ref{rmk:infinite}.
\end{Rmk}

\begin{Rmk}
	Another important preliminary remark is that the eigenvalues of a smooth $(1,1)$-tensor (or a $(0,2)$-tensor) $T$ depend continuously on the base point on a $m$-manifold. Indeed, this is a standard result coming from the continuity of roots of a polynomial. In particular if at a point we have signature $(a,b)$ that is $a$ positive eigenvalues, $b$ negative ones and $m-a-b$ zeros, then in a neighborhood of the point, we have \textit{at least} $a$ positive eigenvalues and $b$ negative ones (the amount of zeroes might have decreased). This idea will be key for the lemmas of this section.
\end{Rmk}

Now here is the last lemma.

\begin{Lem}\label{lem:fin}
	\hspace{\parindent} Let $\Omega$ be a Kobayashi hyperbolic domain in $\Mc^n_\kappa,n>3$, and $\Sigma$ a hypersurface with non-vanishing second fundamental at a point $p\in\Sigma$. Then there exists a totally geodesic immersion $\pi:\Mc_\kappa^3\rightarrow\Mc_\kappa^n$ passing through $p$ and $\pi^{-1}(\Sigma)$ has non-zero mean curvature at $p$.
\end{Lem}

\begin{proof}
	Let us start by a quick remark on curvature in a transversal intersection. Assume at a point $q\in\Sigma$ we choose two principal directions $e_1$, $e_2$ with principal curvatures $\kappa_1$ and $\kappa_2$ respectively and let $\nu(q)$ be the normal of $\Sigma$ at $q$.
	
	The key property we will use is that for every $3$-dimensional subspace of $T_p\Mc_\kappa^n$ for a given point $p$, there exists a unique totally geodesic immersion $\iota:\Mc_\kappa^3\rightarrow\Mc_\kappa^n$ tangent to it. Indeed, the existence is straighforward:
	\begin{itemize}
		\item[--] In $\RR^n$, one can simply take the corresponding affine $3$-plane.
		\item[--] In $\SS^n$ (and $\HH^n$), one can take the intersection of $\SS^n$ (respectively $\HH^n$) with the affine $4$-plane in the ambient space determined by the given $3$-space and the normal direction.
	\end{itemize}
	The uniqueness also follows easily from the symmetry of the model spaces, but we will not need this fact here.
	
	Therefore let $\pi:\Mc_\kappa^3\rightarrow\Mc_\kappa^n$ be the totally geodesic immersion passing through $q$ tangentially to $\Span(e_1,e_2,\nu(q))$. The intersection between $\Sigma$ and $\pi(\Mc_\kappa^3)\cap\Omega$, near $q$, is tranverse. In particular $\pi^{-1}(\Sigma)$ is a surface. The important remark here is that the two principal curvatures of $\pi^{-1}(\Sigma)$ at $\pi^{-1}(q)$ are still $\kappa_1$ and $\kappa_2$. This is easy to show because the normal in $\pi^{-1}(\Omega)$ of $\pi^{-1}(\Sigma)$ at $\pi^{-1}(q)$ is still $\nu(q)$ by definition of $\pi$. Therefore, the second fundamental form in $\pi^{-1}(\Omega)$, at $q$, it is just the restriction of the second fundamental form in $\Omega$. Note that this is false at other points because the unit normal in $\pi^{-1}(\Sigma)$ has no reason to stay the same as the one in $\Omega$. To summarize, the key idea is that we can take a totally geodesic slice and reduce the situation to dimension $3$ while preserving the principal curvatures at the point of interest.
	
	With this observation, the remaining part of the proof becomes easy. Firstly, if there are two principal curvatures of same sign at $p$, then the proof is over by taking both eigendirections, and the normal one for a slice. Secondly, if there are not two principal curvatures of the same sign, then we are in one of the three situations:
	\begin{itemize}
		\item[--] Either there is one positive curvature and $n-2$ zeroes at $p$, or
		\item[--] there is one negative curvature and $n-2$ zeroes at $p$, or
		\item[--] there are one positive curvature, one negative and $n-3$ zeroes at $p$.
	\end{itemize}
	To conclude, we can just take one non-zero eigendirection with one zero. Note that $n>3$ is required for the third situation so that a zero eigendirection exists. In any case, the resulting surface will have non-vanishing mean curvature at $p$.
\end{proof}

\begin{Rmk}
	We could rephrase Lemma \ref{lem:ex} in the following way: a non-minimal surface in a Kobayashi hyperbolic domain in $\Mc_\kappa^3$ is at finite distance. Moreover Lemma \ref{lem:fin} could also be rephrased as: any non-totally geodesic hypersurface in $\Mc_\kappa^n,n>3$ can be sliced into a non-minimal surface in $\Mc_\kappa^3$. This last sentence obviously fails with minimal surface already in $\Mc_\kappa^3$.
	
	Moreover, recall that the strong maximum principle states that given two minimal hypersurfaces in a Riemannian manifold which are tangent at a point and staying on one side from each other must coincide. In a way both lemmas show the failure of this principle when one of the submanifolds is not a hypersurface. However, in dimension three, a minimal disk is a hypersurface. In particular, the strong maximum principle shows that it is impossible to find a strictly separated minimal disk to a minimal hypersurface in a $3$-manifold.
	
	Therefore, the only thing we can say in $\Mc_\kappa^3$ is that a non-minimal surface is at finite distance. This remark highlights the intriguing behaviour appearing in dimension $3$.
\end{Rmk}

We are now ready to prove Theorem \ref{thm:2}. The proof is immediate from Lemma \ref{lem:ex} and Lemma \ref{lem:fin}.

\begin{Thm*}
\hspace{\parindent} Let $(\Mc,g)$ be a Kobayashi hyperbolic Riemannian manifold with constant sectional curvature and dimension $n>3$, $\Sigma$ a submanifold and $p\in\Sigma$. If $\Sigma$ is not, near $p$, a totally geodesic hypersurface then there is a point $p'\in\Sigma$ near $p$ which is at finite Kobayashi distance.
\end{Thm*}

\begin{proof}
	Let us start by assuming that $\Sigma$ is a hypersurface and that $\Mc=\Omega$ a Kobayashi hyperbolic domain of a model space. Note that by hypothesis, in any neighborhood in $\Sigma$ of $p$ there exists a point $p'$ with non-vanishing second fundamental form. In particular, $\Sigma$ follows the conditions of Lemma \ref{lem:fin} at $p'$ therefore there exists a totally geodesic immersion $\pi:\Mc_\kappa^3\rightarrow\Mc_\kappa^n$ passing through $p'$ such that $H_{\pi^{-1}(\Sigma)}(p')\neq0$.
	By Lemma \ref{lem:ex}, there is a sequence of point $a_k$ in $\pi^{-1}(\Omega)$ converging to $\pi^{-1}(p')$ and a fixed point $b\in\pi^{-1}(\Omega)$ such that there is a constant $C$ with $d_{\pi^{-1}(\Omega-\Sigma)}(a_k,b)\leq C$ for all indices $k$. In particular, taking $q:=\pi(b)$, $p_k:=\pi(a_k)$ and using the fact that $\pi$ is a totally geodesic immersion, we have that $d_{\Omega-\Sigma}(p_k,q)\leq d_{\pi(\Mc_\kappa^3)\cap(\Omega-\Sigma)}(p_k,q)\leq d_{\pi^{-1}(\Omega-\Sigma)}(a_k,b)\leq C$. Having $p_k\rightarrow p'$, $p'$ is at finite distance.
	Now, we just need to study the case where $\Sigma$ is not a hypersurface. This is done re-using the same ideas in the following proposition.
\end{proof}

\begin{Pro}\label{pro:last}
	\hspace{\parindent} Let $\Omega$ be a Kobayashi hyperbolic domain in $\Mc^n_\kappa,n\geq 3$, and $\Sigma$ be a submanifold such that $\codim\Sigma>1$. Then every point of $\Sigma$ is at finite distance.
\end{Pro}

\begin{proof}
	Let us start by reducing the situation similarly to the other lemmas.
	
	Firstly, assume $\Sigma$ is not a point and let $p\in\Sigma$, $\nu_1$, $\nu_2$ be two linearly independant transversal vectors at $p$ and $e$ be a non-zero tangent vector to $\Sigma$. Let $\pi:\Mc_\kappa^3\rightarrow\Mc_\kappa^n$ be the slice passing through $p$ tangentially to $\Span(\nu_1,\nu_2,e)$. By construction, $\pi^{-1}(\Sigma)$ is a smooth regular curve.
	
	Secondly, if $\Sigma$ is a point, we cannot take $e$ tangent. However, we can take $e$ to be a third linearly independant transversal vector and the resulting submanifold after the slice will be a point.
	
	Note that it is enough to construct a strictly separated sequence of disks when $\pi^{-1}(\Sigma)$ is a curve. Indeed, if it was a point, any strictly separated sequence of disks for a curve containing $\pi^{-1}(\Sigma)$ will also be strictly separated for the point. Therefore it is enough to show Proposition \ref{pro:last} when $n=3$ and $\Sigma$ a smooth regular curve.
	
	Because $\Sigma$ is regular, there exists an infinite amount of surfaces (in a neighborhood of $p$) which contains $\Sigma$. Let $S$ be one of them. If $H_S(p)=0$ then we can slightly deform $S$ while preserving $\Sigma$ to have $H_S(p)\neq 0$. In particular, Lemma \ref{lem:ex} shows that $p$ is at finite distance from $\Omega-S$. By the decreasing property of the Kobayashi-Royden pseudometric: $p$ is at finite distance from $\Omega-\Sigma$.
\end{proof}

\begin{Ex}
Let us finally give a class of examples of minimal hypersurfaces at finite distance to illustrate the previous proofs. Let $S$ be a minimal surface (non totally geodesic) in $\RR^3$, $\Sigma_n:=S\times\RR^{n-3}$ is a minimal hypersurface in $\Omega$ a bounded domain of $\RR^n,n>3$ (assuming that $\Omega\cap\Sigma_n\neq\varnothing$). Then by Theorem \ref{thm:2}, $\Sigma^n$ is at finite distance because it is non totally geodesic. These examples illustrate well the proof of both lemmas. Indeed at a point $p\in\Sigma_n$ we can take the slice $p+\Span(e_1,e_3,e_n)$. The resulting surface is cylindrical, that is the product of strictly convex line and $\RR$ in $\RR^3$. In this situation, it is easy to find a separated catenoid that can be put on top of the surface.
\end{Ex}

\section*{Conclusion}

\hspace{\parindent} We now have a better understanding of what type of submanifolds works well with Kobayashi's notion of hyperbolicity. But there are still a lot of open questions and reflexions to have relative to those questions. In particular:
\begin{itemize}
	\item[--] Is it possible to extend Theorem \ref{thm:2} on an arbitrary (Kobayashi hyperbolic) Riemannian manifold ? Are the totally geodesic hypersurfaces the only submanifolds at infinite distance ?
	\item[--] What happens in dimension $3$ ? Especially even the simple example of a catenoid in a bounded Euclidean domain is still open. Moreover do we have: \emph{$\Sigma$ at infinite distance if and only if it is totally geodesic} or \emph{$\Sigma$ at infinite distance if and only if it is minimal} ?
	\item[--] Could we obtain an analogous Green's theorem ? See \cite{gr72} and \cite{gr77}. Knowing now what kind of submanifolds are detected by the Kobayashi-Royden pseudometric, a Green's theorem would be the natural next step. Specifically, how many totally geodesic hypersurfaces do we need to remove on a compact Riemannian manifold for it to become complete hyperbolic ?
\end{itemize}

\bibliographystyle{siam}
\bibliography{ref}
	
\end{document}